\documentclass{amsart}

\usepackage{mathtools,enumerate}
\usepackage[square,compress,comma, numbers,sort]{natbib}
\usepackage{color}
\usepackage{mathrsfs}   
\usepackage[shortlabels]{enumitem}
\usepackage{cleveref}
\usepackage[T1]{fontenc}

\newcommand{\R}{\mathbb{R}}

\newcommand{\eps}{\varepsilon}

\newcommand\bkny[1]{ \begin{eqnarray*} #1 \end{eqnarray*}}

\newcommand\ve{\varepsilon}

\newcommand{\abs}[1]{\left\lvert #1 \right\rvert}

\DeclarePairedDelimiterXPP\pk[1]{\mathbb{P}}\{ \}{}{ #1}
\DeclarePairedDelimiterXPP\E[1]{\mathbb{E}}\{ \}{}{	#1}

\DeclarePairedDelimiterXPP\ind[1]{\mathbb{I}}( ){}{	#1}

\newcommand\toprob{\overset{p}\rightarrow}

\newcommand\cE[1]{#1}
\definecolor{c30}{rgb}{0.,0.,1.}
\def\cEE#1{#1}

\numberwithin{equation}{section}
\newtheorem{theorem}{Theorem}[section]
\newtheorem{lem}[theorem]{Lemma}
\newtheorem{korr}[theorem]{Corollary}
\theoremstyle{remark}
\newtheorem{remark}[theorem]{Remark}

\numberwithin{equation}{section}

\newcommand{\prooftheo}[1]{ \textsc{Proof of Theorem} \ref{#1} }

\newcommand{\prooflem}[1]{\textsc{Proof of Lemma} \ref{#1}}
\newcommand{\proofkorr}[1]{\textsc{Proof of Corollary} \ref{#1}}

\newcommand{\QED}{\hfill $\Box$}

\newcommand{\COM}[1]{}

\newcommand\IF{\infty}

\newcommand{\inr}{\in \R}

\newcommand{\BQN}{\begin{eqnarray}}
\newcommand{\EQN}{\end{eqnarray}}
\newcommand{\BQNY}{\begin{eqnarray*}}
\newcommand{\EQNY}{\end{eqnarray*}}
\newcommand\ldot{, \ldots,}

\newcommand{\limit}[1]{\lim_{#1 \to   \infty}}

\newcommand\todis{\overset{d}\rightarrow}

\newcommand\bqny[1]{ \begin{eqnarray*} #1 \end{eqnarray*}}
\newcommand\bqn[1]{ \begin{eqnarray} #1 \end{eqnarray}}

\newcommand{\BS}{\begin{sat}}
\newcommand{\ES}{\end{sat}}
\newcommand{\BT}{\begin{theorem}}
\newcommand{\ET}{\end{theorem}}
\newcommand{\BK}{\begin{korr}}
\newcommand{\EK}{\end{korr}}

\newcommand{\BEX}{\begin{example}}
\newcommand{\EEX}{\end{example}}

\newcommand{\BD}{\begin{de}}
\newcommand{\ED}{\end{de}}

\newcommand{\BIT}{\begin{itemize}}
\newcommand{\EIT}{\end{itemize}}
\newcommand{\BDI}{\begin{description}}
\newcommand{\EDI}{\end{description}}

\newcommand{\BRM}{\begin{remark}}
\newcommand{\ERM}{\end{remark}}

\newcommand{\BEL}{\begin{lem}}
\newcommand{\EEL}{\end{lem}}

\newcommand{\nelem}[1]{{Lemma \ref{#1}}}

\newcommand{\netheo}[1]{{Theorem \ref{#1}}}

\newcommand{\nekorr}[1]{{Corollary \ref{#1}}}

 \newcommand\Z{\mathbb{Z}}
\newcommand\inn{\in \mathbb{N}}

\newcommand{\set}[1]{\left\{#1\right\}}

\newcommand\MAF{\mathcal{M}(\setA,X)}

\newcommand\setA{\mathscr{A}}
\newcommand\MA{\mathcal{M}(\setA,X)}

\newcommand\Iz{\mathcal{J}_z}
\newcommand\LFA{\mathcal{L}(\setA,f(X))}

\newcommand\Lz{\mathcal{L}(\setA, \Iz(X))}

\newcommand\LH{\mathcal{L}(\setA,f_{\kappa}(X))} 
\newcommand\setD{D}
\newcommand\TTT{\R^d}
\newcommand\TT{\R^d}

\newcommand\pka{p_\kappa}

\newcommand\LHH{\mathcal{L}(\setA, I_z)} 
\newcommand\LHz{\mathcal{L}(\setA,  \Iz)}

\begin{document}

\title{The harmonic mean formula for random  processes}

\author{Krzysztof Bisewski}
\address{Department of Actuarial Science, University of Lausanne, Switzerland}
\email{Krzysztof.Bisewski@unil.ch}
\author{Enkelejd Hashorva}
\email{Enkelejd.Hashorva@unil.ch}
\author{Georgiy Shevchenko}
\address{Kyiv School of Economics, Ukraine}
\email{gshevchenko@kse.org.ua}


\begin{abstract}
	Motivated by the harmonic mean formula in  \cite{MR1009865}, 
	we investigate  the relation between the sojourn time and supremum of a random process $X(t),t\inr^d$ and extend the harmonic mean formula for general  stochastically continuous $X$. We  discuss two applications concerning the continuity of distribution of supremum of $X$  and representations of classical  Pickands constants.
\end{abstract}

\maketitle

\section{Introduction}
Let $X(t),t\in \TTT$ be a measurable and separable real-valued random process.  Given a measurable function $f: \R \mapsto\R$, define the $f$-sojourn time (or $f$-occupation time) of $X$ in the set $\setA\subset \TTT$  by 
$$
\LFA= \int_\setA f(X(t)) \lambda(dt),
$$
where $\lambda$ is the Lebesgue measure on $\TTT$ and further $A$ is $\lambda$-measurable. When $\setA$ is countable, then we simply take $\lambda$ to be the counting measure on $\setA$.
{An example of $f$, which is particularly interesting in application} is 
 $$f(t)=\cEE{\Iz}(t):= \ind{t> z},\ z\inr,$$
where $\ind{\cdot}$ is the indicator function. \cEE{Hereafter we set $\MA= \sup_{t\in \setA } X(t)$. 
Note in passing that  when $\lambda(A)$ is finite, both  $\MA$ and $\Lz$ afe determined by finite-dimensional distributions (fidi's) of $X$ (see  \cite[Lem 10.4.2]{MR3561100}). } 
 Consider \cEE{next} the simpler case that  $\setA$ is a countable subset of $\R^d$. 
Clearly, for all $z\inr$ satisfying 
    $\pk{\MA>z}>0$ we have 
 \begin{equation}\label{eq0} 
 \pk{\MA > z, \Lz=0}=0.	
\end{equation} 
An   implication of  \eqref{eq0} is an elegant \cEE{and important} result of Aldous, the so called harmonic mean formula, which is
stated   in \eqref{eq1} below. 
\BEL[\protect{\cite[Lem 1.5, Eq. (1.7)]{MR1009865}}] Let $\setA=\{t_1 \ldot t_n\} \subset  \TTT, n>1$ be given and let $X(t),t\in \TTT$ be a real-valued random process.  If $z\inr$ is such that $\pk{\MA>z}>0$,  then  
\begin{equation}\label{eq1}
	 \pk{\MA > z} = \sum_{i=1}^n \pk{X(t_i) > z } \E*{ \frac{1}{\Lz} \Bigl \lvert X(t_i)> z  }.
\end{equation}
\label{lem0}
\EEL 

For general $\setA$ and $X$  the claim in \eqref{eq0} does not hold in general. 
 For instance, considering a deterministic process
$X(t)=0,t\in [0,1]\setminus \{ 1/2\}$ and $X(1/2)=1$ we clearly have that 
$$\sup_{t\in [0,1]} X(t)=\sup_{t\in (0,1)} X(t)=1>0$$
 almost surely but 
$\pk{\int_0^1 X(t)\lambda(dt) =0}=1$.

Two natural \cEE{questions} which arise here are:
\begin{enumerate} [Q1)]
	\item \label{Qa}Under what   conditions \cE{on fidi's} of $X$ does \eqref{eq0}   hold?
	\item \label{Qb} \cE{Can  the mean harmonic formula be extended }to more general $f$ and uncountable $\setA$?
\end{enumerate}

 \nelem{lem0} is of interest in several applications. In the present form or with some adjustments, it has been utilised extensively for  rare events simulations (by importance sampling technique), see, e.g., \cite{ADB,MR3844167}.
We shall show in the next section that a natural assumption under which \Cref{Qa} has a positive answer is that of stochastically continuous $X$,
\cE{which is  satisfied if, for instance, $X$ is separable, jointly measurable and has stationary increments, see Remark \ref{GennaPartha}, \Cref{item4:GennaPartha} below.} Moreover, for stochastically continuous $X$ we shall  derive the harmonic mean formula in \netheo{thG}. \\
The manuscript is organized as follows. Our main findings are presented in Section \ref{s:main_result}, where we also discuss two applications concerning the continuity of the distribution of $\MA$ and derive new representations for the classical Pickands constants. All the proofs are relegated to Section \ref{s:proofs}.

\section{Main Result}\label{s:main_result}
In this section we shall consider $X(t), t\in \TTT$ stochastically continuous, i.e.\ for all $t\in \TTT$ we have the convergence in probability 
$X(s) \toprob X(t)$ as $s\to t$. In view of \cite[Thm 1, p.\ 171]{MR636254} this guarantees that there exists a jointly measurable and separable version. Therefore in the following, whenever $X$ is stochastically continuous we shall further assume  that $X$ is jointly measurable and separable, see \cite[Thm 1]{MR356206}, \cite[Thm 9.4.2]{MR1280932} for equivalent conditions that guarantee measurability of a random process.
Below we shall consider  open $\setA\subset   \TTT$ 
answering  both questions raised in the Introduction. 
\begin{theorem}
Let $f:\R \mapsto \R$ be  a \cEE{measurable function}  and let $z\in \R, \kappa \in [-\IF, z]$ be given. Suppose that  
 $\setA\subset \TTT$ 
	is  open  and $X(t),t\in \TTT$ is a  stochastically continuous real-valued process. 
	If $f(x)> 0$ for all 
	$x>\kappa$ \cE{and $\pk{\MAF>z}>0$}, then 
\begin{equation}\label{eth0}
	\pk{ \MAF >z , \LH =0} =0,
\end{equation}	
where $f_{\kappa}(x)= f(x)\ind{x >\kappa}, \ x\in \R$. 
If further $\LH$ is almost surely finite, then  
	\bqn{\label{eth1}
		\pk*{ \MAF > z} 
		&=&    \int_{\setA}  \E*{ \frac{ f_\kappa(X(t))\ind{\MAF> z} }  
			{\LH } } \lambda(dt).
} 
\label{thG}	
\end{theorem}

A simple application of \netheo{thG} concerns the question of continuity of $\MAF$, which has been investigated in various generalities in numerous contributions, see e.g., \cite{tsirel1976density,MR1604537,LifBook,MR4046509} and the excellent contribution  \cite{MR1604537}, 
where the methodology is explained in  details.

\BK Let $   X(t), t\in \R^d$ be  as in \netheo{thG}, let $z\inr$ be given and suppose the open set $\setA\subset \TTT$ is bounded. If \cEE{$\pk{\MAF \ge z}>0$}, $\pk{X(t)=z}=0$ for all $\cEE{t\in \setA}$, and 
  \bqn{ \label{lifsh} 
  	\pk{\MAF \ge z, \LHH=0}=0, 
 }
\cEE{where $I_z(t):=\ind{t\ge z}$,} then the distribution of $\MAF$ is continuous at $z$. 
\label{thG2}	
\EK

\BRM \begin{enumerate}[(i)]
	\item Theorem~\ref{thG} is also true for countable $\setA$, with $\lambda$ being the counting measure.
Moreover, 	the assumption that $\setA$ is open can be relaxed to $\setA$ is a Borel subset of $\R^d$ with non-empty 
  interior  $\setA^o$  satisfying  $\pk{\sup_{t\in \setA} X(t) = \sup_{t\in \setA^{o}} X(t)}=1$; 
  \item\label{item:drop_indicator} If $\kappa=z$, then the indicator function in \eqref{eth1} can be dropped; 
	\item For $X(t),t\in \TTT$ measurable and separable with separant $\setD$, the assumption of stochastic continuity can be relaxed to the following: For all  $\ve>0$ positive and given $t_i \in \setD$  
\cE{\bqn{	\label{demain}
	\pk{ X(t) \le \cE{z}, X(t_i)> z+\cE{\ve}} \to 0,  \   t \to t_i. 
	}
}
\cE{To see that, note that in the proof of  \netheo{thG} below  the crucial argument of the proof (along with stochastic continuity) is
	\eqref{crucial}.}  
In particular, \eqref{demain} holds if for all $\eps>0$
	$$
	\pk{X(t)<X(t_i)-\eps}\to 0, \ t\to t_i,
	$$
	or  $X$ is stochastically continuous from the right when $d=1$. For $d>1$ the latter assumption can be formulated in terms of quadrant stochastic continuity;
	\item \cE{If   $X$ is Gaussian with continuous covariance function on $\TTT$, then  \eqref{demain} is satisfied;}  	
	\item Two common choices for  $f$ in the results above are $f(x)= e^{bx}$ and $f(x)= \abs{x-z}^b, b\ge 0$, for which we have that $f(x)>0$ for all $x>z$. In particular, from 
	\eqref{eq0} for all $b\ge 0$  we have	
	\begin{equation}\label{kol}
		\pk*{ \sup_{t\in \setA} \abs{X(t)}> 0, \int_\setA \abs{X(t)}^b \lambda(dt)=0}=0.
	\end{equation} 
	\item \label{item4:GennaPartha} When $X$ is stationary with $\pka=\pk{X(0)> \kappa}>0$, then $\pka=\pk{X(t)> \kappa}, \forall t\in \TTT$, and thus, under the assumptions of \netheo{thG}, by the shift-invariance of Lebesgue measure  we can rewrite \eqref{eth1} as  
	\bqn{\label{eth2}
		\pk*{ \MAF > z} 	&=&   \pka \int_{\setA}   \E*{ \frac{ f(X(0)) \ind{\max_{s\in \setA} X(s-t)> z}}  {\int_\setA f(X(s-t)) \ind{ X(s-t)> \kappa} \lambda(ds) }   \Bigl\lvert  X(0)> \kappa } \lambda(dt).
	} 
Moreover, if we assume that $X$ is separable and jointly measurable, then by \cite[Prop 3.1]{Roy1} (see  \cite{MR3561100}[Thm 1.3.3] for the more general case of processes with stationary increments) we have that $X$ is stochastically continuous.\\
	If $\setA \subset  \delta \Z^d, \delta >0$ and $\lambda$ is the counting measure on $\delta \Z^d$ we have that \eqref{eth2} still holds;
	\item If we take $X(t)=2\ind{U< 0}+ \sin(t-U) \ind{ U\ge 0}, t\in [-1,1]$, with $U$ uniformly distributed on $[-1,1]$, it follows that for any open $\setA \subset [0,1]$
 $$
  \pk*{ \MAF \ge  1, \int_0^1 \ind{ X(t) \ge 1 } dt > 0 }=1/2, 
  $$
  hence  condition \eqref{lifsh} in \nekorr{thG2} cannot be removed. 
  This example was kindly suggested by Mikhail Lifshits (personal communication).
 \end{enumerate}
\label{GennaPartha}
\ERM

Next, we present another application of the harmonic mean formula.  The idea, which appears for instance in \cite{Berman82,Berman92}  is to find some positive functions $v(x), q(x),x >0$ 
such that 
$$Y_\kappa(t)= v(\kappa)[X(q(\kappa)t) - \kappa], \ t\in \TTT,$$
conditioned on  $Y_\kappa(0)>0$ converges weakly to some random process $Y(t),t\in \TTT$, as $\kappa \to \IF$ .\\
In the sequel, let $X(t),t\ge 0$ be a centered stationary Gaussian process with unit variance  and correlation function $r$ 
satisfying the so-called Pickands condition see \cite{Pit96}, i.e., 
\begin{equation}\label{Pic}
 r(t)< 1, \forall t \ge 0, \ 1- r(t)= C \abs{t}^\alpha (1+ o(1)), \ t\to 0
\end{equation}
 for some $C>0, \alpha \in (0,2]$. In view of \nelem{lemA} below  for this case we can take  
 $$v(\kappa)=\kappa, \ q(\kappa)= \kappa^{-2/\alpha}$$
  and the limit process $Y$ is given with $B$ \cEE{being} a standard fractional Brownian motion with Hurst index $\alpha/2$  by 
\begin{equation}\label{gbY}
	Y(t)= W_1(C^{1/\alpha}t)+ \eta, \quad  W_b(s)= b[ \sqrt{2} B(t)- \abs{t}^{\alpha}], \quad b\ge 0,t\inr ,
\end{equation} 	
 where  the unit exponential rv  $\eta$ is independent of $B$.

In the following, for $\delta>0$ let $\delta\Z$ denote the infinite grid of uniformly spaced points $\delta\Z := \{\ldots,-2\delta,-\delta,0,\delta,2\delta,\ldots\}$ and set $0\Z = \R$.
\BT If $X(t),t\ge 0$ is a centered stationary Gaussian process with unit variance, \cEE{continuous sample paths}  and correlation function satisfying the Pickands condition for some $C>0, \alpha \in (0,2]$, and $T_z>0$ \cEE{is such that} $\limit{z}T_z z^{2/\alpha}=\IF$, then for all any $\delta\geq0$
\bqn{\label{electri} 
	\pk*{ \sup_{t\in [0,T_z]\cap\delta\Z} X(t)> z} \sim C^{1/\alpha} \mathcal H^\delta_\alpha  T_z  z^{2/\alpha} \pk{X(0)> z}  , \quad z\to \IF, 
}
provided that the left-hand side above converges to 0 as $z\to \IF$,
with
\bqn{\label{lastF} \mathcal H_\alpha^\delta = e^{\theta}\E*{ \frac{ \ind{  \sup_{t\in \delta \Z} W_1( t)+ \eta > \theta} }{  \delta\sum_{t\in \delta \Z} e^{   W_b(  t)}
			\ind{ W_1( t)+ \eta > 0}     }},
}
for all $\delta>0$, and
\bqn{
\label{PickA:2}
\mathcal H_\alpha := \mathcal H^0_\alpha = \lim_{\delta\to0} \mathcal H^\delta_\alpha = e^{\theta} \E*{ \frac{\ind{  \sup_{t\inr} W_1( t)+ \eta > \theta} }{  \int_\R e^{   W_b(  t)}\ind{ W_1( t)+ \eta > 0}\lambda(dt) } } ,
}
where $b, \theta \in [0,\IF)$  are arbitrary.
\label{satE}
\ET

\BRM \begin{enumerate}[(i)]
	\item 
 The new results above are the  expressions for Pickands constant $H_\alpha$,  which are known for the case $a=  \theta=0,$ see
  \cite[3.6]{HBernulli} and \cite[Thm 1.1]{MR4029237}, see also \cite[J20a,~J20b]{MR969362};
 \item The discrete Pickands constant is defined for $X$ as in \Cref{satE} by  
\bqn{ \label{haD} 
	\mathcal H_\alpha^\delta= 
	\limit{T} T^{-1} \E*{\sup_{ t\in [0,T] \cap \delta \Z} e^{ W_1(t)} } = 	\limit{T}  \limit{z} 
	\frac{ \pk{  \sup_{t\in [0,T] \cap \delta \Z} X( q(z) t ) > z}}{ T \pk{X(0)> z} },
}
with $q(z) \sim z^{-2/\alpha}$ as $z\to \IF$. \cE{See e.g. \cite[p.\ 1605]{DiekerY} for the first formula and 
 \cite[p.\ 164]{MR2154366} for the second \cEE{one}.}
Furthermore, we have  
\bqn{\label{lastF} \ \ \  \delta \mathcal H_\alpha^\delta &=&  e^\theta\E*{ \frac{ \ind{  \sup_{t\in  \delta \Z} W_1(  t)+ \eta > \theta} }{ \sum_{t\in 
 			\delta \Z} e^{   W_b(   t)}
 		\ind{ W_1( t)+ \eta > 0}     }}	=  \E*{ \frac{   \sup_{t\in  \delta \Z} e^{ W_1(  t)} }{ \sum_{t\in 
 				\delta \Z} e^{   W_1(   t)}}} \\
 			&	=  &\pk*{ \sup_{ t  \in \delta \mathbb{N}  } W_1(  t ) + \eta \le 0 } >0
 			\label{lastF2}
}
 	for all $\delta >0$. The first formula for $\theta=0$  {and both two other \cEE{ones} for  $\mathcal H_\alpha^\delta$} in \eqref{lastF} and    \eqref{lastF2} can be derived also as in \cite{MR3745388} utilising previous findings  of \cite{BojanS}, see also \cite{Hrovje, BP,PH2020,kulik:soulier:2020, Planic}. The expression in  \eqref{lastF2} was derived in \cite{Albin1990} and appeared also latter in \cite{MR1772400,AlbinC};
 	\item 
 	\cE{Explicit formulas for $\mathcal H_1^\delta, \delta>0$ are known, see \cite[Lem 5.16, Rem 5.17]{KW}. Such a formula appeared also in other connections, see 
 	definition of $\nu(x)$ function in \cite[Eq.~(2.1)]{MR1331667}. An alternative \cEE{expression}  is given in \cite{buritic2021variations}.} 
 	\end{enumerate}
\ERM

\section{Proofs}\label{s:proofs}

\prooftheo{thG} 
We show first that 
\bqn{ \label{isra}
	\pk*{ \MAF >0, 	 \int_{\setA} \ind{X(t)>0}\, \lambda(dt) =0 }=0.
}
Suppose for simplicity that $d=1$ and consider without loss of generality  $X$ to be further measurable and separable (we use definition in \cite{MR636254} for (joint) measurability).  Let $\setD \subset\R $ be a countable dense set that  is a separant for $X$   and set $\setD\cap \setA = \set{t_i,i\inn }$.  
For fixed arbitrary $\eps>0$, by definition of separability we have  
$$\MAF=\sup_{t\in \setA \cup D} X(t) = \sup_{i \inn } X(t_i)$$
almost surely. 
By the  stochastic continuity of $X$, for each positive integer $i$ and given $\ve>0$
$$
\pk*{|X(t) - X(t_i)|> \eps }\to 0,\quad t\to t_i.
$$
Therefore, since $\setA$ is open,  there exists  some open interval $\Delta_i\subset \setA$ containing $t_i$   such that 
\bqn{ \label{crucial} 
	\pk*{ {X(t)\le 0}, X(t_i)>\eps } \le \eps 2^{-i}, \ \ \forall t\in \Delta_i.
}	  
Consequently, for all $i\inn$, the Fubini-Tonelli theorem yields
$$
\E*{\int_{\Delta_i} \ind{X(t)\le 0}\lambda(dt)\ind{X(t_i)>\eps}} = \int_{\Delta_i} \pk*{ X(t)\le 0, X(t_i)>\eps }\lambda(dt)  \le \eps 2^{-i}\lambda(\Delta_i).
$$
Hence, by  the Markov inequality
\begin{align*}
	\pk*{ {\int_{\Delta_i} \ind{X(t)>0}\,dt = 0},X(t_i)>\eps }& = \pk*{ {\int_{\Delta_i} \ind{X(t)\le 0}\lambda(dt) = \lambda(\Delta_i)}, 
		X(t_i)>\eps } \\
	&\le \pk*{ {\int_{\Delta_i} \ind{X(t)\le 0}\lambda(dt) \ge  \lambda(\Delta_i)},X(t_i)>\eps }\\
	& \le \frac{1}{\lambda(\Delta_i)} \E*{\int_{\Delta_i} \ind{X(t)\le 0}\lambda(dt)\ind{X(t_i)>\eps}} \\
	& \le \eps 2^{-i}. 
\end{align*}
Since further 
\bkny{ 
  		\pk*{\MAF >\eps, \int_{\setA} \ind{X(t)>0}\lambda(dt) = 0}
	&=& \pk*{ {\sup_{i\ge 1} X(t_i)>\eps} , {\int_{\setA} \ind{X(t)>0}\lambda(dt) = 0}}\\
	&\le &
	\sum_{i\ge 1} \pk*{X(t_i)>\eps , {\int_{\Delta_i} \ind{X(t)>0}\,dt = 0}}\\
	&\le &\eps \sum_{i\ge 1}2^{-i} = \eps
}	
the claim  in \eqref{isra} follows by letting $\eps\to 0$. \\
Next, by our assumption, $f_\kappa(X(t))=f(X(t)) \ind{X(t)> \kappa}>0$  whenever $X(t)>\kappa$, thus the equality 
$$\LH=\int_\setA f(X(t)) \ind{X(t)> \kappa} \lambda(dt)=0$$ \cEE{implies} 
 $$\lambda(\set{t\in \setA: X(t) >\kappa})=0.$$
  Therefore, applying \eqref{isra}  to the process $X(t) - \kappa, t\in \TTT$ 
we obtain
\bkny{
	\lefteqn{
	 	\pk{ \MAF >z , \LH =0}}\\
 	&\le &\pk{ \MAF >z ,\lambda(\set{t\in \setA: X(t) >\kappa}) =0}\\
	& \le& \pk*{ \MAF >z , \int_{\setA}\ind{X(t)>\kappa}\lambda(dt) = 0}\\
	& \le &\pk*{ \MAF >z , \int_{\setA}\ind{X(t)>z}\lambda(dt) = 0} \\
			& = & 0,
}
 which establishes the first claim.

Next, by the assumption, 
$$\LH= \int_{\setA} f(X(t)) \ind{ X(t)>\kappa} \lambda(dt) $$
is almost surely finite and non-zero, interpreting $\frac 0 0$  as 0, \cE{by \eqref{eq0}} we find that
$$  \mathbb{I}(\MAF> z )= 	
\frac{\LH} {\LH } \mathbb{I}(\MAF > z )	
$$
 holds  almost surely. Hence by the Fubini-Tonelli theorem and the assumption $\pk{\MAF>z}>0$    
\bqny{
	\pk{ \MAF > z}  &=&\E*{ \frac{\LH}{\LH} \ind{\MA  >z }	  	}
	\\
	&=&\int_{\setA}   \E*{  \frac{ f( X(t))  \ind{\cEE{\MA} > z} \ind{X(t)> \kappa} } { \LH }} \lambda(dt)  
	>0
}
establishing   the proof. 
\QED

\proofkorr{thG2}.
\cEE{Recall that $\Iz(s)= \ind{s> z}$ and we set $I_z(s)=\ind{s\ge z}$.} Since we suppose that  
\bqn{\label{eI} 
	\pk{ \ind{ X(t)=z}=0}=1, \quad \forall t\in \setA
}
and $\setA$ is bounded,  then 
$$\IF>  \LHH=\int_\setA\ind{ X(t)\ge  z }\lambda(dt) \ge \LHz=  \int_\setA\ind{ X(t)> z }\lambda(dt)$$
almost surely and hence 
$$ \E{\LHH-\LHz }= \int_\setA \E{   \ind{ X(t)=z}}\lambda(dt) =0$$
implying 
$$\pk{\LHH=\LHz}=1.$$
Since, by the assumption, $\MAF\ge z$ implies almost surely  $\LHz\in (0,  \IF)$,  as in the proof of \netheo{thG}, \cEE{using further \eqref{eI}} we have
\bqny{ 
	\pk*{ \MAF \ge  z} 
	=    \int_{\setA}  \E*{ \frac{ \ind{X(t) \ge z}  }  
		{\LHH } } \lambda(dt)= \int_{\setA}  \E*{ \frac{ \ind{X(t)> z}  }  
		{\LHz } } \lambda(dt)=\pk*{ \MAF > z} ,  
} 
where the last equality follows  {from Theorem \ref{thG} and Remark \ref{GennaPartha}, \Cref{item:drop_indicator}, This concludes the proof.}
\QED 
\\

Next we state a lemma which is needed in the proof of \netheo{satE}.  Denote by $C([0,T])$ the space of real-valued 
continuous functions on $[0,T]$ equipped with a metric which turns it into a Polish space and let $\mathcal{C}$ be the corresponding Borel $\sigma$-field. 
We write $\todis$ for the weak convergence of fidi's. Hereafter 
$W_b(t)$ and $Y(t)=W_1(t)+ \eta$ are as in  \eqref{gbY},  with $\eta$ a unit exponential rv independent of $W_1$. Set n the following 
\begin{equation}\label{def:H_T}
\mathcal H^{\delta}_\alpha(T) := \lim_{z\to\infty} \frac{\pk{\sup_{t\in\delta\Z\cap[0,T]}X(tz^{-2/\alpha}) > z}}{\pk{X(0)>z}} 
\end{equation}
\cEE{for some $T>0$.}
\BEL Let the stationary process $X(t),t\in \cEE{[-T,T],T>0}$ have almost surely sample paths in $C([-T,T])$. Suppose further that  $X$ is  Gaussian with mean zero and unit variance function satisfying   \eqref{Pic} for some $C>0$ 
and $\alpha\in (0,2]$ and let $q(z)$ be positive such that $q(z) \sim z^{-2/\alpha}$ as $z\to \IF$.  
\begin{enumerate}[(i)]
	\item\label{lemA:A}  As $z \to \IF$ we have the weak convergence $
z[ X(q(z)t)- z] \mid \big( X(0)> z\big)	\todis   Y(t)$;
	\item \label{lemA:B} For all $b,\theta \ge 0$ and all $\delta, T$ positive 
\begin{equation}\label{eq:to_show_truncated_pickands}
\mathcal H^{\delta}_\alpha(T) = {e^\theta} \sum_{\tau\in\delta\Z\cap[0,T]} \E*{\frac{\ind{  \sup_{  t\in\delta\Z\cap[0,T]} W_1( t-\tau) + \eta > \theta } }{ \sum_{t\in\delta\Z\cap[0,T]} e^{ W_b(t-\tau)} \ind{ W_1(t-\tau) + \eta > 0}}}.
\end{equation}
\end{enumerate}
\label{lemA}
\EEL 

\prooflem{lemA} \Cref{lemA:A}: The convergence of fidi's is well-known, see for instance \cite[Lem 2]{AlbinC}. \\
\Cref{lemA:B}: For any fixed $\theta>0$ set $\kappa(z) := z-\theta/z, z>0$. In the following, for brevity we write $\kappa = \kappa(z)$. Take 
\begin{equation}\label{eq:def_f_kappa}
f_\kappa(x)= e^{b \kappa x }\ind{x>\kappa}, \quad q_\theta(x) := q(\tfrac{x+\sqrt{x^2+4\theta}}{2})
\end{equation}
such that $q_\theta(\kappa(z)) = q(z)$ and $Y_\kappa(t)=\kappa(z)[X(q_\theta\big(\kappa(z)t\big)- \kappa(z)]$. 
Using Theorem \ref{thG}, with $f_\kappa$ defined in \eqref{eq:def_f_kappa} we obtain
\begin{align*}
\mathcal H^{\delta}_\alpha(T) = \lim_{z\to\infty} \frac{\sum_{\tau\in\delta\Z\cap[0,T]}\E{U_z(\tau,T)}}{\pk{X(0)>z}}, \quad U_z(\tau,T) := \frac{f_\kappa(X(q(z)\tau))}{\sum_{t\in K}f_\kappa(X(q(z) t))}\ind*{ \sup_{s\in K} X( q(z) s)> z}.
\end{align*}
For any $\tau\in\delta\Z\cap[0,T]$ we have
\begin{align*}
U_z(\tau,T) & = \frac{e^{b\kappa X(q_\theta(\kappa)\tau)} \ind*{ X(q_\theta(\kappa)\tau) > \kappa}}{\sum_{t\in \delta\Z\cap[0,T]}e^{b\kappa X(q_\theta(\kappa) t)} \ind*{ X(q_\theta(\kappa)t) > \kappa}}\ind*{ \sup_{t\in \delta\Z\cap[0,T]} \kappa(X( q_\theta(\kappa) t)-\kappa) > \kappa(z-\kappa)} \\
& = \frac{e^{ b\kappa(X(q_\theta(\kappa)\tau)-\kappa)} \ind*{ \kappa(X(q_\theta(\kappa)\tau)-\kappa) > 0}}{\sum_{t\in \delta\Z\cap[0,T]}e^{ b\kappa(X(q_\theta(\kappa)t)-\kappa)} \ind*{ \kappa(X(q_\theta(\kappa)t)-\kappa) > 0}}\ind*{ \sup_{t\in \delta\Z\cap[0,T]} \kappa(X( q_\theta(\kappa) t)-\kappa)>\theta\kappa/z}.
\end{align*}
Further by the the stationarity of $X$ 
\begin{align*}
U_z(\tau,T) \overset{d}{=} \frac{e^{bY_\kappa(0)}\ind*{ Y_\kappa(0) > 0}}{\sum_{t\in \delta\Z\cap[0,T]}e^{bY_\kappa(t-\tau)}  \ind*{ Y_\kappa(t-\tau) > 0}}\ind*{ \sup_{t\in \delta\Z\cap[0,T]} Y_\kappa(t-\tau)>\theta} := U^*_z(\tau,T).
\end{align*}
\cEE{Hence}  applying  \Cref{lemA:A} and the continuous mapping theorem we obtain the convergence in distribution
\begin{align}\label{eq:def_U}
\left(U^*_z(\tau,T) \mid Y_\kappa(0) > 0\right)& \overset{d}{\to} \frac{\ind{  \sup_{  t\in\delta\Z\cap[0,T]} W_1( t-\tau) + \eta > \theta } }{ \sum_{t\in\delta\Z\cap[0,T]} e^{ W_b(t-\tau)} \ind{ W_1(t-\tau) + \eta > 0}}:= U(\tau,T), \ z\to \IF.
\end{align}
Moreover, since $U^*_z(\tau) \leq 1$, then by the dominated convergence theorem $\lim{z}\E{U_z(\tau,T)}= \E{U(\tau,T)}.$ Consequently, using further 
$$\pk{Y_\kappa(0) > 0} = \pk{X(0)> z - \theta/z}\sim  e^\theta \pk{ X(0)> z}, \quad z\to \IF$$ 
implies 
\begin{align*}
\mathcal H^{\delta}_\alpha(T) = \lim_{z\to\infty} \frac{\pk{Y_\kappa(0)>z}}{\pk{X(0)>z}}\sum_{\tau\in\delta\Z\cap[0,T]}\E{U^*_z(\tau)\mid Y_\kappa(0)>0} = e^\theta\sum_{\tau\in\delta\Z\cap[0,T]}\E{U(\tau,T)}
\end{align*}
establishing proof. 
\QED

\prooftheo{satE}. The claim in \eqref{electri} is shown \cEE{in \cite[Lem 9.15]{Pit19}}, \cite[Thm 2.1]{MR4127347} for $\delta=0$. The proof  for $\delta>0$ follows with similar arguments and is omitted. Alternatively, 
 it can be derived utilising the approach in \cite[p.~339]{AlbinC} and applying  \nelem{lemA}. 
In view of   
\cite[Thm 1.1]{MR4029237} and \cite[Eq.~5.2]{HBernulli} we have the following expressions 
$$ \mathcal  H_\alpha=  \E*{ \frac{\ind{  \sup_{t\inr} W_1( t)+ \eta >0} }{ \int_\R \ind{ W_1( t)+ \eta > 0}\lambda(dt) }     } = 
  \lim_{ \delta \downarrow 0} \mathcal  H_\alpha^\delta =\lim_{ \delta \downarrow 0} \frac{1}{\delta} \E*{ \frac{ \ind{  \sup_{t\in \delta \Z} W_1( t)+ \eta > 0} }{ \sum_{t\in \delta \Z} 
		\ind{ W_1( t)+ \eta > 0}     }}>0.   $$
Note that the positivity of the constants follows from the fact that  almost surely (see e.g., \cite[Thm~6.1, Prop~6.1]{WangStoev})
\bqn{ \label{eqG}
	W_1(t) \to - \IF,\  t\to \IF.
}
Further, $\mathcal H_\alpha = \lim_{\delta\downarrow 0} \mathcal  H_\alpha^\delta$ is well-known, see e.g., \cite{ZKE}. 
Consequently,  \eqref{PickA:2} follows for $b=0$ and all $\delta,\theta\ge 0.$ 
Next,   in view of \cite{Htilt}[Thm 2.8] or \cite[Eq.~(5.4)]{MartinE} (see also \cite{PH2020} for case $d=1$)  for all $b, \theta \in [0,\IF), x>0 $ we have (write $B^h f(\cdot) = f(\cdot-h), h\inr$) 
 \bqn{ 
 	 \label{eqR12_prim} \E{ F( x e^{ \eta + W_1 } ) \ind{   W_1(h)+ \eta  >- \ln x}} = 
	x \E{  F(   e^{\eta+B^h W_1}) \ind{ W_1(-h)+ \eta  >\ln x }   }, \forall h\in \TT,
}
valid for all measurable maps $F:C([0,\IF))\mapsto [0,\IF]$; 
\cEE{we interpret $0/0$ and $0 \cdot \IF$ as 0}.  
\cEE{For some countable dense set $D \subset \R$ set}
$$L(t)= e^{W_1(t)+ \eta}, \ M= \cEE{\sup_{t\in \R \cap D}} L(t) ,\  S= \int_{\R} L(t)\lambda(dt), \quad 
\chi=e^\theta.
$$
\cEE{Observe that}
$$L(s)/M\in (0,1], \ M \in (0,\IF), \ S \in (0,\IF), \  \int_\R  [L(t)]^b  \ind{ L(t)>1 }  dt \in (0,\IF)$$
almost surely (the last claim follows from  \eqref{eqG}, the previous last follows from  \cite[Thm~6.1, Prop~6.1]{WangStoev}).
\cEE{Hence, using \eqref{eqR12_prim},    we obtain from the Fubini-Tonelli theorem that}
\bkny{
\lefteqn{\E*{ \frac{ e^\theta \ind{\sup_{t\in \R } W_1(t)+ \eta > \theta} }{ \int_\R e^{W_b(t)}\ind{ W_1( t)+ \eta > 0}\lambda(dt) }     }}\\
	 &=& \chi \E*{ \frac{S}{S} \frac{ [L(0)]^ b \ind{ \sup_{\cEE{t\inr \cap D}} L(t) >\chi }} {  \int_\R  [L(t)]^b  \ind{ L(t)>1 }  \lambda(dt)  }}\\
	& = & \chi \int_\R \E*{ \frac{ [L(0)]^ b    L(s)  \ind{ M >\chi }   } { S \int_\R   [L(t)]^b  \ind{ L(t)>1 } \lambda(dt)  }  } \lambda(ds)\\
	& = & \chi \int_0^\IF \int_\R \E*{ \frac{ [L(0)]^ b   \ind{ rL(s)> 1}  \ind{ M >\chi }    } { S \int_\R   [L(t)]^b  \ind{ L(t)>1 } \lambda(dt) } } r^{-2} \lambda(ds) dr.
}
Write $M^*=rM$  and $S^*=rS$,  $L^*= r L$. Then the expectation above can be written as
$$\E*{ \frac{r [L^*(0)]^b \ind{ r L(s)> 1} \ind{  M^*> r \chi}}{ S^*\int_\R [ L^* (t)]^b  \ind{L^*(t)> r} \lambda(dt)}}  =
r \E*{ G(r L)   \ind{ r L(s)> 1}},$$
where
$$G(f) := \frac{f^b(0)\ind{\sup_{\cEE{t\in\R \cap D}}f(t)>r\chi}}{\int_\R f(t)dt \int_\R f^b(t)\ind{f(t)>r}\lambda(dt)}.$$
After applying \eqref{eqR12_prim} for any $r>0$ we obtain
$$r\E*{ G(r L)   \cEE{\ind{ r L(s)> 1}}} = r^{2} \E*{ G( B^s L) \cEE{\ind{L(-s)>r}}},$$
which, after applying substitution $z := r\chi/M$, further gives us
	\bkny{
	\lefteqn{\chi \int_0^\infty\int_\R \E*{ G( B^s L) 1(L(-s)>r)}\lambda(ds) dr}\\
	&=&   \chi  \int_\R \int_0^\IF \E*{ \frac{ [L(-s)]^ b   \ind{  L(-s)> r}  \ind{ M > r\chi }    } { S\int_\R   [L(t)]^b  \ind{ L(t)>r } \lambda(dt) }  } \lambda(ds) dr\\
	&=&     \int_0^\IF\int_\R   \E*{ \frac{ M [L(-s)]^ b   \ind{   \chi  L(-s)/M> z , 1> z }    } { S\int_\R   [L(t)]^b  \ind{  \chi  L(t)/M>z } \lambda(dt) }  } 
	\lambda(ds) dz\\
	&=&   \int_0^1   \E*{ 		\frac{ M } {S} \cdot \frac{\int_\R  [L(-s)]^ b \ind{\chi   L(-s)/M>z } ds}{ \int_\R   [L(t)]^b  \ind{ \chi L(t)/M>  z } \lambda(dt) }    
}  dz\\
		&=& \E*{ \frac{ \sup_{ t\in \R} e^{W_1(t)} } { \int_\R e^{W_1(t)}   \lambda(dt)}} =  \mathcal  H_\alpha.
}
In the last equality we used that $z< 1$ and $L(s)/M< 1$ almost surely implies $z< L(s)/M \chi$ almost surely for all $s\inr$ since $\chi\ge 1$, hence the proof for the case $\delta=0$ in complete.
\cEE{If $\delta>0$ the above calculations can be repeated. An alternative proof follows by passing to the limit in the expression given in \eqref{eq:to_show_truncated_pickands}.}
 \QED  

\section*{Acknowledgments}
\cEE{We thank the reviewer for several important suggestions.}

\bibliographystyle{ieeetr}
\bibliography{EEEA}

\end{document}